\documentclass[12pt]{amsart}

\usepackage{amsmath,bm}
\usepackage{amssymb}
\usepackage{graphicx}
\usepackage{upgreek}  
\usepackage{mathrsfs} 
\usepackage{color}
  
\newtheorem{theorem}{Theorem}[section]

\newtheorem{lemma}[theorem]{Lemma}
\newtheorem{proposition}[theorem]{Proposition}
\theoremstyle{definition}

\newtheorem{remark}[theorem]{Remark}

\newcommand{\eq}[1]{(\ref{#1})}
\newcommand{\wt}[1]{\widetilde{#1}}
\newcommand{\wh}[1]{\widehat{#1}}

\newcommand{\linesunder}[3]{\LSU{\begin{array}[t]{c}\underbrace{#1}\vspace*{.5em}\end{array}}{\mbox{\footnotesize\rm #2}}{\mbox{\footnotesize\rm #3}}}

\newcommand{\LSU}[3]{\begin{array}[t]{c}#1\vspace*{-1em}\\_{#2}\vspace*{-.5em}\\_{#3}\end{array}}

\def\vv{{\bm v}}
\def\uu{{\bm u}}
\def\bb{{\bm b}}
\def\hh{{\bm h}}
\def\nn{{\bm n}}

\def\EE{{\bm e}}
\def\vvt{\vv_{\text{\sc t}}^{}}
\def\wtvvt{\wt{\vv}_{\text{\sc t}}^{}}
\def\divS{{\rm div}_\text{\sc s}^{}}
\def\nablaS{\nabla_\text{\sc s}^{}}
\newcommand{\R}{\mathbb R}
\newcommand\DT[1]{\mathchoice
                 {{\buildrel{\hspace*{.1em}\text{\LARGE.}}\over{#1}}}
                 {{\buildrel{\hspace*{.1em}\text{\Large.}}\over{#1}}}
                 {{\buildrel{\hspace*{.1em}\text{\large.}}\over{#1}}}
                 {{\buildrel{\hspace*{.1em}\text{\large.}}\over{#1}}}}
\newcommand{\ITEM}[2]{\parbox[t]{.055\textwidth}{{\rm #1}}\hfill\parbox[t]{.945\textwidth}{#2}\vspace*{.8mm}}

\newcommand{\TTT}{\color{black}}
\newcommand{\FFF}{\color{black}}
\newcommand{\SSS}{\color{black}}
\newcommand{\EEE}{\color{black}}

\numberwithin{equation}{section}

\title[Optimal control of semi-compressible fluids]{Optimal control of flows of
  viscoelastic semi-compressible fluids}

\author[T. Roub\'\i\v cek]{Tom\'a\v s Roub\'\i\v cek}

\address{\hspace*{-1.2em}Mathematical Institute, Charles University,\newline
Sokolovsk\'a 83, CZ-186 75 Praha 8
\newline
{\rm and}
\newline
Institute of Thermomechanics\\Czech Academy of Sciences,\newline
Dolej\v skova 5, CZ-182 00 Praha 8, Czech Republic}

\email{{\tt tomas.roubicek@mff.cuni.cz}}



\begin{document}
\begin{marginpar}

\allowdisplaybreaks

\begin{abstract}
  Semilinear parabolic systems with bi-linear
  nonlinearities cover a lot of applications and their
  optimal control leads to relatively simple optimality
  conditions. An example is the incompressible Navier-Stokes system
  for homogeneous fluids, which is however here
  modified towards a physically reasonable model
  of slightly (so-called ``semi'') compressible 
  liquids rather than fully compressible gases. An optimal control problem optimizing
  also pressure on the boundary 
  is considered and, in the simple variant, analysed as far as uniqueness
  of the control-to-state mapping and 1st-order optimality conditions
  in the 2-dimensional case and outlined in a nonsimple variant for the
  3-dimensional case. Some other bi-linear parabolic systems as Cahn-Hilliard
  diffusion or magneto-hydrodynamics can be treated analogously.
\end{abstract}

\maketitle


\section{Introduction}
  The optimal control of semilinear parabolic systems with bi-linear
  nonlinearities allows to exploit such special nonlinear structure,
  in particular to formulate optimality conditions in a relatively
  lucid way. An example of such bi-linear nonlinearities arises
  from convective time derivatives in Eulerian description of
  continuum-mechanical models, as properly used in fluid dynamics.
  The bi-linear nonlinear terms thus occurs in incompressible Newtonien
  homogeneous fluid flows, described by the incompressible Navier-Stokes
  equations.

Optimal control of such system has been vastly addressed in
literature mainly in the fully incompressible variants.
Standardly expected results are available especially in   
two-dimensions, cf.\ \cite{AbeTem90SCPF,FatSri94NSCO,HinKun01SOMO}.
In three dimensions, rather only steady-state situations have been scrutinized
(as in \cite{MalRou99OSFI,ReyTro07OCSN,RouTro03LSOC,TroWac06SOSO})
while only particular results are available in the evolutionary situations
\cite{AbeTem93OCPG,Casa93OCPG,FatSri92EOCV,KiRoWa17PPOC} because of the
well-recognized difficulties related with lack of regularity and
uniqueness of the weak solutions.

Incompressible fluids are however only an idealized model
and, although well applicable in many situations, it ignores various
physical phenomena (most importantly the propagation of pressure waves)
and, as mentioned, brings even serious analytical difficulties
especially in three-dimensional cases. On the other hand, the
fully compressible models involves nonlinearities of more complicated
structure than only bi-linear. Actually, although most fluids (liquids or
melted metals \TTT or magma\EEE) are quite compressible, they are not really ``fully''
compressible as gases. E.g.\ water is about 50-times more compressible
than steel but, nevertheless, it is far from to be so compressible as gases
and has a specific density even under zero pressure \TTT and even can withstand
a certain negative pressure (cf.\ e.g.\ \cite{CADA12EWOL} and references therein)\EEE,
in contrast to gases. 

To reflect this phenomenology, models for a class of so-called
{\it semi-compres\-sible fluids}
have been devised in \cite{Roub20QISI} as a compromise between fully
compressible models with substantially varying mass density 
and the fully incompressible Navier-Stokes model.
They pursue the following attributes:
\begin{enumerate}
\item[$(\upalpha)$] 
propagation of \TTT longitudinal waves (i.e.\ pressure waves, called also P-waves) \EEE 
  is allowed and their dispersion is controlled in a certain way,
\item[$(\upbeta)$] the energy balance is preserved at least formally,
  but in some models even rigorously,
\item[$(\upgamma)$] the pressure is well defined
  also on the boundary,
\item[$(\updelta)$] the equations are consistently written in Eulerian
  coordinates (i.e.\ the model is fully convective),
\item[$(\upepsilon)$] in some models, uniqueness of weak solutions holds even in the physically
  relevant 3-dimensional cases.
\end{enumerate}
For optimal control, it is important that such models still exhibit bi-linear
nonlinear structure. In particular, we will avoid usage of continuity equation
for mass density\TTT. 
This \EEE is well acceptable for fluids whose density is
(nearly) constant in space and varies (nearly) negligibly with pressure\FFF, so that
it can be modeled as constant\EEE.  
Let us remind that the involvement \TTT of \EEE density as a variable
governed by the continuity equation $\DT\varrho+{\rm div}(\varrho\vv)=0$ would
make \TTT the \EEE term $\varrho(\vv\cdot\!\nabla)\vv$ in \eq{NS1+} below
tri-linear, which would make optimality conditions more complicated\TTT, not
mentioning many other analytical problems including uniqueness of the
state response\EEE.

In Section~\ref{sec-system}, we first specify a simple semi-compressible model
and an initial-boundary-value problem for it.
Some analytical properties of
it are then analyses in Section~\ref{sec-n=2} for the two-dimensional case
and then in Section~\ref{sec-control} we consider a simple optimal-control
problem with a quadratic cost and formulate first-order optimality conditions. 
Eventually, in Section~\ref{sec-remarks}, we outline some modifications
towards three-dimensional case and enhancements of the basic model to obtain
some better properties or applicability to other phenomena coupled with the
fluid flow\TTT{s}\EEE.

\section{Semi-compressible fluids}\label{sec-system}
\def\nablax{\nabla}
\def\nablaxx{\nabla^2}
\def\pp{p}
\def\GM{\gamma}
\def\CHI{{\bm\vartheta}}
\def\CHIi{{\bm\vartheta}_i}
\def\CHIj{{\bm\vartheta}_j}
\def\LAMBDA{\pi}
\def\faa{f}
\def\fbb{g}
\def\II{I}
\def\G{\varGamma}
\def\O{\varOmega}
\def\d{{\rm d}}
\def\In{\in}
\def\Rn{\mathbb R^n}
\def\e{\varepsilon}
\def\qqqquad{\qquad\qquad}

We consider a fixed bounded domain $\varOmega\subset\Rn$ with a
Lipschitz boundary $\varGamma$ and a finite time interval $I=[0,T]$.
We will use the dot-notation means the partial derivative in time.
The basic scenario will use the velocity $\vv$ and pressure $\pp$, and consider
a model devised in \cite{Roub20QISI}. \TTT It consists of the system of two
parabolic equations epressing momentum equation and convective transport and diffusion of pressure\EEE:
\begin{subequations}\label{NS-eqn-evol+}\begin{align}\label{NS1+}
&\varrho\DT\vv
+\varrho(\vv\cdot\!\nablax)\vv-{\rm div}(\nu\EE(\vv))
+\frac\varrho2({\rm div}\,\vv)\,\vv+
\nablax\Big(\pp+\frac\beta2\pp^2\Big)
\,=\,\uu
\!\!\!
&&\text{on }\II{\times}\O\,,
\\[-.2em]\label{NS2+}
&\beta
(\DT\pp+\vv\cdot\!\nablax\pp)
+{\rm div}\,\vv=\GM\Delta\pp&&\text{on }\II{\times}\O\,,
\intertext{\TTT completed by the boundary and initial conditions\EEE}
&
\big[\nu\EE(\vv)\nn\big]_{\text{\sc t}}^{}\!+b\vvt
=\,0\,,\ \ \  
\nn{\cdot}\vv=\,0,\ \
\text{ and }\ \ \nn{\cdot}\nablax\pp
=0
&&\text{on }\II{\times}\G\,,
\label{NS3+}
\\[-.2em]\label{NS4+}
&\vv(0,\cdot)=\vv_0\ \ \ \text{ and }\ \ \ \pp(0,\cdot)=\pp_0&&\text{on }\O\,
\end{align}\end{subequations}
with \FFF mass density $\varrho>0$ assumed constant (in particular independent of
pressure), with \EEE some constants $\beta$ and $\GM>0$ \TTT commented below\EEE,
viscosity $\nu>0$, and $\EE(\vv)=\frac12\nabla\vv
+\frac12(\nabla\vv)^\top$. The bulk force $\uu$ is here prescribed and
later in Section~\ref{sec-control} will be used as a distributed control.
Let us note that,
beside the
usual ``hydrostatic''  pressure $\pp$, there is also the pressure contribution
$\frac12\beta\pp^2$ in \eq{NS1+} due to the elastic internal energy of
the fluid. The system \eq{NS-eqn-evol+} gets a good physical consistency in Eulerian
description with $\beta>0$ being the impressibility. In terms of the
bulk elastic modulus $K$ in physical units Pa=J/m$^3$, the impressibility is
$\beta=1/K$. This modulus determines the velocity of 
\TTT P\EEE-waves
(namely $\sqrt{K/\varrho}\,$=\,sound speed, provided $\GM=0$) which can
propagate through such fluids, in contrast to ideally incompressible fluids.
Moreover, $\GM>0$ allows for modeling a
certain dispersion of pressure waves and is motivated by a mass
diffusion in the continuity equation in the full compressible model,
advocated by H.\,Brenner \cite{Bren05KVT,Bren06FMR}. The physical
dimension of $\GM/\beta$ is m$^2$/s and, vaguely speaking,
dividing it by a ``characteristic velocity'' of the flow
and a ``characteristic size'' of the system,
one gets a dimensionless P\'eclet number (or, in fluid dynamics,
also called Brenner's number) expressing dominance of either
the convective or the diffusive transport phenomena.
\TTT The mentioned compromise between fully
incompressible models (where P-waves cannot propagate at all) and fully
compressible models (where density can substantially vary)
is well legitimate when pressure variations (and thus density variations in
compressible situations) are much smaller than the elastic bulk modulus $K$
and simultaneously the shear elastic modulus is zero, which is relevant in most
situation in liquids (but usually not in gases neither in solids). 
E.g.\ water has the elastic bulk modulus
about 2\,GPa (which is much less than e.g.\ rocks with $>$10\,GPa or steel with
$>$100\,GPa) but still it is much larger than
usual pressure variations in most practical situations, and the speed of P-waves
is about 1.5\,km/s. These numbers are known with 5-digit (or more) accuracy
(including its dependence on temperature, salinity, and pressure itself,
cf.\ IAPWS-standard \cite{WagPru02IAPW})
and putting $K=\infty$ (which is, in fact, done in
incompressible models) is a simplification which might be often not well
acceptable.\EEE

The (so-called Navier) boundary conditions \eq{NS3+} involves the tangential
velocity $\vvt=\vv-(\nn{\cdot}\vv)\nn$.

This semi-compressible modification now contributes also to the stored
energy due to the calculus
\begin{align}\label{quasi-incompress-energy}
\int_\varGamma\pp(\vv{\cdot}\nn)\,\d S&-\int_\varOmega\nablax\pp\cdot \vv\,\d x
=\int_\varOmega\pp\,{\rm div}\,\vv\,\d x
\\&\nonumber=\int_\varOmega\pp\,\Big(\beta\DT\pp
+\beta \vv\cdot\!\nablax p-\GM\Delta\pp
\Big)\,\d x
\\&=\int_\varOmega\GM|\nablax\pp|^2-\frac{\beta}2\pp^2{\rm div}\,\vv\,\d x
+\frac{\d}{\d t}\int_\varOmega
\frac{\beta}2\pp^2\,\d x
\nonumber\end{align}
if employing also the boundary condition $\nn{\cdot}\vv=0$
and $\nablax\pp\cdot\nn=0$. More specifically, for the term
$\beta\pp \vv\cdot\!\nabla\pp$, we have used the Green formula for
\begin{align*}
 \int_\varOmega\pp \vv\cdot\!\nabla\pp\,\d x
&=\int_\varGamma\pp^2\vv{\cdot}\nn\,\d S-\int_\varOmega{\rm div}(\pp \vv)\pp\,\d x
  \\[-.0em]&=\int_\varGamma\pp^2\vv{\cdot}\nn\,\d S
  -\int_\varOmega(\nabla\pp\cdot \vv)\pp+\pp^2{\rm div}\,\vv\,\d x
  \\&=\int_\varGamma\frac12\pp^2\vv{\cdot}\nn\,\d S
  -\int_\varOmega\frac12\pp^2{\rm div}\,\vv\,\d x\,.
\end{align*}
This storage-energy mechanism together with the kinetic
energy just facilitates wave propagation.
The extra bulk force $\frac\varrho2({\rm div}\,\vv)\,\vv$ in \eq{NS1+},
proposed by R.\,Temam \cite{Tema69ASEN},
arises by (slight) compressibility and is presumably small
as ${\rm div}\vv$ is presumably very small but (slightly) violates Galilean
invariance of the model, as pointed out in \cite{Toma??ITST}.
This is the price payed for omitting the continuity equation for $\varrho$
and simplifying considerably the analysis of the model, which also gets
bi-linear semi-linear structure with still keeping a lot of physically
relevant features. This extra ``structural'' force
balances the energetics due to the calculus
\begin{align}\label{convective-tested}
\int_\varOmega\varrho(\wh \vv\cdot\!\nabla)\vv{\cdot} \vv+
\frac12\int_\varOmega\varrho|\vv|^2({\rm div}\,\wh \vv)\,\d x
=\frac12\int_\varGamma\varrho|\vv|^2(\wh \vv{\cdot}\nn)\,\d S\,
\end{align}
to be used for $\wh \vv=\vv$. The overall energy balance looks as
\begin{align}\label{energy-NS-semi}
&\int_\varOmega\!\!\!\!\linesunder{\frac\varrho2|\vv(t)|^2_{_{_{_{_{_{_{}}}}}}}}{kinetic}{energy}\!\!\!\!\!\!+\!\!\!\!\!\!\linesunder{\frac{\beta}2\pp(t)^2_{_{_{_{_{_{_{}}}}}}}}{elastic}{energy}\!\!\!\!\!\d x
 +\int_0^t\!\int_\varOmega\!\!\!\linesunder{\nu|\EE(\vv)|^2+\GM|\nablax\pp|^2_{_{_{_{_{_{_{}}}}}}}}{dissipation rate}{in the bulk}\!\!\!\d x\d t
  \\[-.1em]&\nonumber\ \ \ \ \ +\int_0^t\!\int_\varGamma
  \!\!\!\!\!\!\!\!\!\!\!\!\linesunder{
  b|\vvt|^2_{_{_{_{_{_{_{}}}}}}}}{dissipation rate}{on the boundary}\!\!\!\!\!\!\!\!\!\!\!\!
  \d S\d t 
=\int_0^t\!\int_\varOmega\!\!\!\!\!\!\!\!\linesunder{
  \uu{\cdot}\vv_{_{_{_{_{_{_{}}}}}}}}{power
  of}{the control}\!\!\!\!\!\!\!\!\!\d x\d t
  +\int_\varOmega\!\!\!\!\linesunder{\frac\varrho2|\vv_0|^2+\frac{\beta}2\pp_0^2}{\TTT initial kinetic and\EEE}{\TTT stored energy\EEE}\!\!\!\!\!\,\d x\,.
\nonumber\end{align}
In fact, \eq{energy-NS-semi} holds for weak solutions rigorously only for $n=2$
while in higher dimensions, it might hold only as an inequality unless the
solution is enough regular.

\begin{remark}[Quasi-incompressible fluids]
Without the convective term $\beta \vv\cdot\!\nablax\pp$ in \eqref{NS2+}
and the corresponding pressure contribution $\frac12\beta\pp^2$
in \eq{NS1+}, the system (\ref{NS-eqn-evol+}a,b)
models so-called quasi-incompressible fluids, considered
as an artificial regularization of the incompressible
model rather for numerical purposes, cf.\ e.g.\ R.\,Temam
\cite[Ch.\,III,\,Sect\,8]{Tema77NSE} or, with $\GM>0$,
A.\,Prohl \cite{Proh97PQCM},
which can also be understood as a singularly perturbed variant
(mathematically understood as a regularization) of the incompressible
model, i.e.\ $\beta=0$ and $\gamma=0$. For analytical purposes, the
quasi-incompressible regularization with $\beta>0$ and $\gamma>0$ was
devised by A.P.\,Oskolkov \cite{Osko73SPQL}.
\end{remark}

\section{The state problem in two-dimensions}\label{sec-n=2}

Throughout this article, we will use the standard notation for the function
spaces: the Lebesgue and the Sobolev spaces, namely $L^p(\varOmega;\R^n)$ for
Lebesgue measurable functions $\varOmega\to\R^n$ whose Euclidean norm is
integrable with $p$-power, and $W^{k,p}(\varOmega;\R^n)$ for functions from
$L^p(\varOmega;\R^n)$ whose
all derivative up to the order $k$ have their Euclidean norm integrable with
$p$-power. We also write briefly $H^k=W^{k,2}$. Moreover, for a Banach space $X$
and for $I=[0,T]$,
we will use the notation $L^p(I;X)$ for the Bochner space of Bochner
measurable functions $I\to X$ whose norm $\|\cdot\|_X$ is in $L^p(I)$, 
and $H^1(I;X)$ for functions $I\to X$ whose distributional derivative
is in $L^2(I;X)$. Furthermore, $C
(I;X)$ will denote the Banach space of
continuous functions $I\to X$. 
We will use the notation $(\cdot)^*$ for the dual space and define specifically
\begin{subequations}\label{spaces}\begin{align}\label{spaces-V}
    &{\mathcal V}=\big\{\vv\!\in\!L^2(\II;H^1(\O;\Rn))\cap
    H^1(\II;H^1(\O;\Rn)^*);\ \nn{\cdot}\vv\big|_{I{\times}\varGamma}^{}=0
    \big\},\!\!
    \\&{\mathcal P}=\big\{\pp\!\in\!L^2(\II;H^1(\O))\cap
    H^1(\II;H^1(\O)^*)\big\}\,,\ \ \text{ and}
\\&{\mathcal U}=L^2(I{\times}\varOmega;\R^n)\,.
\end{align}\end{subequations}

In any case, under the assumptions 
\begin{align}\label{IC}
\vv_0\In L^2(\O;\Rn)\ \ \text{ and }\ \ \pp_0\In L^2(\O)
\end{align}
with $\uu\in{\mathcal U}$, from \eq{energy-NS-semi} we can read a-priori estimates
$\vv_i,\pp\In L^\infty(\II;L^2(\O))\cap L^2(\II;H^1(\O))$, $i=1,...,n$, for
any dimension. In this and the following section, we will restrict ourselves on $n=2$.

The definition of a weak solution to \eq{NS-eqn-evol+} will be based on the
state-equation mapping $\varPi:
{\mathcal U}\times
{\mathcal V}\times{\mathcal P}\to{\mathcal V}^*\times{\mathcal P}^*$ defined by
\begin{align}\nonumber
\left\langle\varPi(\uu,\vv,\pp),(\wt \vv,\wt\pp)\right\rangle&=
\int_0^T\!\!\!\int_\varOmega\nu\EE(\vv){:}\EE(\wt\vv)
+\Big(\varrho(\vv\cdot\!\nablax)\vv+\frac\varrho2({\rm div}\,\vv)\,\vv
-
\uu\Big){\cdot}\wt\vv
\\[-.2em]&\nonumber
\hspace*{0.2em}+
\Big(\pp{+}\frac{\beta}2\pp^2\Big){\rm div}\,\wt \vv
+\GM\nablax\pp\cdot\!\nablax\wt\pp
+\big(\beta \vv\cdot\!\nablax\pp{+}{\rm div}\,\vv\big)\wt\pp
\,\d x\d t
\\[-.2em]&\nonumber
\hspace*{.2em}+\int_0^T\!\!\!\int_\varGamma 
b\vvt\!\cdot\wtvvt\,\d S\d t
+\int_\varOmega\varrho \vv(T)\cdot\wt \vv(T)+\beta\pp(T)\wt\pp(T)
\\[-.4em]&\hspace*{12.5em}
-\varrho \vv_0\cdot\wt \vv(0)-\beta\pp_0\wt\pp(0)\,\d x\,
\nonumber\end{align}
for any $(\wt\vv,\wt\pp)\in{\mathcal V}^*{\times}{\mathcal P}^*$.
Here we also used that $\mathcal V\subset C(I;L^2(\O;\Rn))$
and $\mathcal P\subset C(I;L^2(\O))$ so that the values
$\vv(t)$, $\wt\vv(t)$, $\pp(t)$, and $\wt\pp(t)$ are well defined in
$L^2$-spaces for $t=T$ or $t=0$. For $\uu\in 
{\mathcal U}$ given, we say that
$(\vv,\pp)\in{\mathcal V}{\times}{\mathcal P}$ is a weak
solution to \eq{NS-eqn-evol+} if $\varPi(\uu,\vv,\pp)=0$.

Existence of weak solutions is by standard arguments: an
approximation e.g.\ by a Galerkin method, usage of apriori estimates to
be read from \eq{energy-NS-semi} written for the approximates solutions,
then passage to the limit by weak convergence and Aubin-Lions compact embedding
theorem for the nonlinear terms.
It is important that, 
if $n=2$,
\begin{subequations}\label{estimates}\begin{align}\label{estimate}
&\|(\vv\cdot\!\nablax)\vv{\cdot}\wt\vv\|_{L^1(I{\times}\O)}^{}\le
\|\vv\|_{L^4(I{\times}\O;\R^2)}^{}\|\nablax\vv\|_{L^2(I{\times}\O;\R^{2\times2})}^{}
\|\wt\vv\|_{L^4(I{\times}\O;\R^2)}^{}
\\&\nonumber
\qqqquad\le
\TTT C_\text{\sc gn}^2\EEE\|\vv\|_{L^\infty(I;L^2(\O;\R^2))}^{1/2}\|\vv\|_{L^2(I;H^1(\O;\R^2))}^{1/2}
\times
\\&\nonumber\qqqquad\quad
\times\|\nablax\vv\|_{L^2(I{\times}\O;\R^{2\times2})}^{}\|\wt\vv\|_{L^\infty(I;L^2(\O;\R^2))}^{1/2}\|\wt\vv\|_{L^2(I;H^1(\O;\R^2))}^{1/2}
\,,
\end{align}
where the Gagliardo-Nirenberg 
\TTT inequality $\|\vv_{12}\|_{L^4(\varOmega;\R^2)}\le
C_\text{\sc gn}\|\vv_{12}\|_{L^2(\varOmega;\R^2)}^{1/2}$ $\|\vv_{12})\|_{H^1(\varOmega;\R^{2})}^{1/2}$ \EEE has been used. Analogously,
one can estimate the term $({\rm div}\,\vv)\vv{\cdot}\wt\vv$.
Also
 \begin{align}\label{estimate+}
   \|\pp^2{\rm div}\,\wt\vv\|_{L^1(I{\times}\O)}^{}&\le
   \|\pp\|_{L^4(I{\times}\O)}^2\|{\rm div}\,\wt\vv\|_{L^2(I{\times}\O)}^{}
   \\&\le \TTT C_\text{\sc gn}^2\EEE
 \|\pp\|_{L^\infty(I;L^2(\O))}^{}\|\pp\|_{L^2(I;H^1(\O))}^{}\|{\rm div}\,\wt\vv\|_{L^2(I{\times}\O)}^{}\,.
\nonumber\end{align}
Similar estimate hold for the nonlinear term $(\vv\cdot\!\nabla\pp)\wt\pp$\,:
\begin{align}\label{estimate++}
&\|(\vv\cdot\!\nablax)\pp{\cdot}\wt\pp\|_{L^1(I{\times}\O)}^{}\le
\|\vv\|_{L^4(I{\times}\O;\R^2)}^{}\|\nablax\pp\|_{L^2(I{\times}\O;\R^{2})}^{}
\|\wt\pp\|_{L^4(I{\times}\O)}^{}
\\&\nonumber
\qqqquad\qqqquad\le
\TTT C_\text{\sc gn}^2\EEE\|\vv\|_{L^\infty(I;L^2(\O;\R^2))}^{1/2}\|\vv\|_{L^2(I;H^1(\O;\R^2))}^{1/2}
\times
\\&\nonumber\qqqquad\qquad\qquad\ \ \
\times\|\nablax\pp\|_{L^2(I{\times}\O;\R^{2})}^{}
\|\wt\pp\|_{L^\infty(I;L^2(\O))}^{1/2}\|\wt\pp\|_{L^2(I;H^1(\O))}^{1/2}.\!
\end{align}\end{subequations}
This shows that indeed $\varPi(\uu,\cdot,\cdot):
{\mathcal V}{\times}{\mathcal P}\to{\mathcal V}^*{\times}{\mathcal P}^*$.

\begin{lemma}[Well-posedness of the controlled system]\label{lem}
  Let $\varrho,\nu,\beta,\gamma>0$, \eq{IC} hold, and $n=2$.
  For any $\uu\in{\mathcal U}$, there is a unique weak solution $(\vv,\pp)\in
  {\mathcal V}{\times}{\mathcal P}$ and the mapping
  $\uu\mapsto(\vv,\pp):{\mathcal U}\to{\mathcal V}{\times}{\mathcal P}$ is
  locally Lipschitz continuous \SSS and also (weak,weak)-continuous\EEE.  
\end{lemma}

\begin{proof}
An important attribute especially in the context of control is uniqueness of
the response $(\vv,\pp)$ for a
given control $\uu$. As for the incompressible Navier-Stokes model, 
the uniqueness holds unfortunately only for two-dimensional
problems. Denoting $\vv_{12}=\vv_1{-}\vv_2$ and $\pp_{12}=\pp_1{-}\pp_2$
for two weak solutions $(\uu_1,\pp_1)$ and $(\uu_2,\pp_2)$
and analysing the identity
$\langle\varPi(\uu,\vv_1,\pp_1)-\varPi(\uu,\vv_2,\pp_2),
(\vv_{12},\pp_{12})\rangle=0$, we have for a.a.\ time instants $t\In I$ (with
$t$ omitted in the following formulas for notational simplicity) that

\begin{align}
\label{unique-est-}
\!\!&
\nu\|\EE(\vv_{12})\|_{L^2(\varOmega;\R^{2\times2})}^2
+\GM\|\nablax\pp_{12}\|_{L^2(\varOmega;\R^2)}^2
+b\|\vv_{12}\|_{L^2(\varGamma;\R^2)}^2
\\\nonumber\!\!&\hspace{12.5em}
+\frac{\d}{\d t}\Big(\frac\varrho2\|\vv_{12}\|_{L^2(\varOmega;\R^2)}^2
+\frac{\beta}2\|\pp_{12}\|_{L^2(\varOmega)}^2\Big)
\\\nonumber\!\!&\hspace{0em}
=\int_\varOmega\!\!
\bigg(\varrho\big((\vv_2{\cdot}\nabla)\vv_2-\vv_1{\cdot}\nabla)\vv_1\big)
\cdot \vv_{12}
+\frac{\beta}2\big(\pp_1^2-\pp_2^2\big){\rm div}\,\vv_{12}
\\[-.5em]\nonumber\!\!&\hspace{1.8em}
+\frac\varrho2\big(({\rm div}\,\vv_2)\vv_2-({\rm div}\,\vv_1)\vv_1\big){\cdot}\vv_{12}
+\beta\big(\vv_1{\cdot}\nabla\pp_1-\vv_2{\cdot}\nabla\pp_2\big)\pp_{12}
\!\bigg)\d x\!\!\!
\\\nonumber\!\!&\hspace{0em}
=\int_\varOmega\!
\!\varrho\Big((\vv_{12}{\cdot}\nabla)\vv_1{+}\frac12({\rm div}\,\vv_1)\,\vv_{12}
{+}(\vv_2{\cdot}\nabla)\vv_{12}
{+}\frac12({\rm div}\,\vv_{12})\,\vv_2
\Big){\cdot}\vv_{12}\,\d x\!\!\!
\\[-.2em]\!\!&\nonumber\hspace{7.em}
+\int_\varOmega\!\beta\Big(\frac{\pp_1{+}\pp_2}2{\rm div}\,\vv_{12}
+\vv_{12}{\cdot}\nabla\pp_1+\vv_2{\cdot}\nabla\pp_{12}
\Big)\pp_{12}\,\d x.\!\!\!
\nonumber\end{align}
The first right-hand integral can be estimated standardly
as for the incompressible Navier-Stokes equation by H\"older's and
Young's inequalities and by the Gagliardo-Nirenberg inequality
combined with Korn's inequality $\|\vv_{12}\|_{L^4(\varOmega;\R^2)}\le
C_\text{\sc gnk}\|\vv_{12}\|_{L^2(\varOmega;\R^2)}^{1/2}\|\EE(\vv_{12})\|_{L^2(\varOmega;\R^{2\times2})}^{1/2}$;
here we again rely on $n=2$.
The same  Gagliardo-Nirenberg inequality holds for $\pp_{12}$ and can be exploited for
estimating the particular terms in the last integral in \eq{unique-est-} as:
\begin{subequations}\label{unique-est}\begin{align}
    &\int_\varOmega
 \frac{\beta}2(\pp_1{+}\pp_2)({\rm div}\,\vv_{12})\pp_{12}\,\d x
\\[-.2em]\nonumber&\quad \le\frac{\beta}2\|\pp_1{+}\pp_2\|_{L^4(\varOmega)}
 \|\vv_{12}\|_{H^1(\varOmega;\R^2)}
 \|\pp_{12}\|_{L^4(\varOmega)}
 \\[-.2em]\nonumber&\quad\le\SSS C_\text{\sc gnk}\EEE
 \frac{\beta}2\|\pp_1{+}\pp_2\|_{L^4(\varOmega)}
 \|\vv_{12}\|_{H^1(\varOmega;\R^2)}
 \|\pp_{12}\|_{L^2(\varOmega)}^{1/2}\|\nablax\pp_{12}\|_{L^2(\varOmega;\R^2)}^{1/2}
\\[-.2em]\nonumber&\quad\le\epsilon \|\vv_{12}\|_{H^1(\varOmega;\R^2)}^2\!
  +\SSS C_\text{\sc gnk}^2\EEE\frac{\beta^2}{16\epsilon}\|\pp_1{+}\pp_2\|_{L^4(\varOmega)}^2
  \|\pp_{12}\|_{L^2(\varOmega)}\|\nablax\pp_{12}\|_{L^2(\varOmega;\R^2)}\!\!
\\[-.2em]&\nonumber\quad\le\epsilon \|\vv_{12}\|_{H^1(\varOmega;\R^2)}^2
+\epsilon\|\nablax\pp_{12}\|_{L^2(\varOmega;\R^2)}^2
 +C\frac{\beta^4}{
   \epsilon^2}\|\pp_1{+}\pp_2\|_{L^4(\varOmega)}^4
  \|\pp_{12}\|_{L^2(\varOmega)}^2\,,
  \\\label{unigue-est-d}&\int_\varOmega\beta
(\vv_{12}{\cdot}\nabla\pp_1)\pp_{12}\,\d x\le\beta
\|\vv_{12}\|_{L^4(\varOmega;\R^2)}\|\nablax\pp_1\|_{L^2(\varOmega;\R^2)}\|\pp_{12}\|_{L^4(\varOmega)}
\\[-.3em]\nonumber&\quad\le
\SSS\beta C_\text{\sc gnk}^2\EEE\|\vv_{12}\|_{L^2(\varOmega;\R^2)}^{1/2}
\|\vv_{12}\|_{H^1(\varOmega;\R^2)}^{1/2}
\times
\\[-.2em]&\nonumber\qqqquad\qqqquad\qquad
\times
\|\nablax\pp_1\|_{L^2(\varOmega;\R^2)}^{}
\|\pp_{12}\|_{L^2(\varOmega)}^{1/2}\|\nablax\pp_{12}\|_{L^2(\varOmega;\R^2)}^{1/2}
\\\nonumber&\quad\le
2\epsilon\|\vv_{12}\|_{H^1(\varOmega;\R^2)}^{}
\|\nablax\pp_{12}\|_{L^2(\varOmega;\R^2)}^{}
\\[-.2em]&\nonumber\qqqquad\qquad\ \
\SSS+C_\text{\sc gnk}^4\TTT\frac{\beta^2}{8\epsilon}\|\nablax\pp_1\|_{L^2(\varOmega;\R^2)}^2
\|\vv_{12}\|_{L^2(\varOmega;\R^2)}^{}\|\pp_{12}\|_{L^2(\varOmega)}^{}\EEE
\\\nonumber&\quad\le
\epsilon\|\vv_{12}\|_{H^1(\varOmega;\R^2)}^2
+\epsilon\|\nablax\pp_{12}\|_{L^2(\varOmega;\R^2)}^2
\\[-.2em]&\qqqquad\quad
+\SSS C_\text{\sc gnk}^4\EEE\frac{\beta^2}{
  \TTT4\EEE\epsilon}\|\nablax\pp_1\|_{L^2(\varOmega;\R^2)}^2\Big(
\|\vv_{12}\|_{L^2(\varOmega;\R^2)}^2+\|\pp_{12}\|_{L^2(\varOmega)}^2\Big)\,,
\nonumber
\\&\int_\varOmega\beta(\vv_2{\cdot}\nabla\pp_{12})\pp_{12}\,\d x\le
  \beta\|\vv_2\|_{L^4(\varOmega;\R^2)}
  \|\nabla\pp_{12}\|_{L^2(\varOmega;\R^2)}\|\pp_{12}\|_{L^4(\varOmega)}
  \\[-.3em]\nonumber&\quad
  \le\beta\SSS C_\text{\sc gnk}^2\EEE\|\vv_2\|_{L^4(\varOmega;\R^2)}
  \|\nabla\pp_{12}\|_{L^2(\varOmega;\R^2)}^{3/2}\|\pp_{12}\|_{L^2(\varOmega)}^{1/2}
  \\&\nonumber\quad\le\epsilon\|\nabla\pp_{12}\|_{L^2(\varOmega;\R^2)}^2
  +
  C\frac{\beta^4}{\epsilon}
  \|\vv_2\|_{L^4(\varOmega;\R^2)}^4\|\pp_{12}\|_{L^2(\varOmega)}^2\,.
\nonumber\end{align}\end{subequations}
Taking $\epsilon>0$ small enough, 
the $\epsilon$-terms on the
right-hand sides of \eq{unique-est} can be absorbed in the left-hand
side of \eq{unique-est-} while the others
can be treated by Gronwall's inequality, using that
$t\mapsto\|\pp_1(t){+}\pp_2(t)\|_{L^4(\varOmega)}^{\FFF4\EEE}$, 
$t\mapsto
\|\nablax\pp_1(t)\|_{L^{\SSS2\EEE}(\varOmega;\R^2)}^{\SSS2\EEE}$
 and $t\mapsto\|\EE(\vv_{\SSS2\EEE}(t))\|_{L^{\SSS4\EEE}(\varOmega;\R^{2\times 2})}^{\SSS4\EEE}$
are $L^1(I)$.

Considering the above two solutions $(\vv_1,\pp_1)$ and $(\vv_2,\pp_2)$ for
two different controls $\uu_1$ and $\uu_2$, respectively, the
right-hand side of \eq{unique-est-} would augment by the term
$\int_\O\uu_{12}\cdot\vv_{12}\,\d x$ with $\uu_{12}=\uu_1{-}\uu_2$.
Then the above estimates in fact show the local Lipschitz continuity of the
control-to-state mapping $\uu\mapsto(\vv,\pp)$ from ${\mathcal U}$
to $(L^2(I;H^1(\O;\Rn))\cap L^\infty(I;L^2(\O;\Rn)))\times
(L^2(I;H^1(\O))\cap L^\infty(I;L^2(\O)))$. By a slight modification
of these estimates, we can also see the local Lipschitz continuity 
$\uu\mapsto(\DT\vv,\DT\pp)$  from ${\mathcal U}$
to $(L^2(I;H^1(\O;\Rn))\cap L^\infty(I;L^2(\O;\Rn)))^*\times
(L^2(I;H^1(\O))\cap L^\infty(I;L^2(\O)))^*$.

\SSS
When $\uu_1$ is fixed and $\uu_2\to\uu_1$ converges weakly in $\mathcal U$, we
can prove the (weak,weak)-continuity by using the compact embedding
of ${\mathcal V}\times{\mathcal P}$ into $L^2(I{\times}\Omega;\R^d{\times}\R)$ by
the Aubin-Lions theorem. This allows to pass to the
limit in all nonlinear terms. At this point, it is also important that $\uu$ occurs
linearly in \eqref{NS1+}.
\end{proof}

Let us end this section by noting that the energy balance \eq{energy-NS-semi}
holds rigorously for any weak solution $(\vv,\pp)$ because $\DT\vv$ is in
duality with $\vv$ and also $\DT\pp$ is in duality with $\pp$, so that
the tests of \eq{NS1+} by $\vv$ and of \eq{NS2+} by $\pp$ are legitimate.

Notably, all these estimates are exact without any ``reserve''.
Let us point out that, in the three-dimensional case,
this semi-compressible model admits only a very weak solution and also 
the uniqueness and continuity of the control-to-state mapping analogous to
Lemma~\ref{lem} is not granted. Cf.\ also Remark~\ref{rem-nonsimple} below for a
modification of the model working for $n=3$.

\section{Optimal control in the two-dimensional case}\label{sec-control}
Beside facilitating pressure waves and their dispersion, the benefit of this
semi-compressible model is that pressure is well defined even with the
traces on the boundary. Also the values of $\pp$ at particular time instants
are well defined in the sense of $L^2(\varOmega)$.
This allows us to involve pressure on the boundary and in the terminal time
into the cost functional which will be considered quadratic for simplicity.
We thus consider the optimal-control problem 
\begin{align}
\!\!\!\left\{\!
\begin{array}{lll}
\mbox{Minimize}\!\!\!\! & 
\displaystyle{\!\!\varPhi(\uu,\vv,\pp):=\!\int_0^T\!\!\!\int_\varOmega
  \frac{\SSS\kappa_1\EEE}2|\vv{-}\vv_{\d}|^2+
  \frac{\SSS\kappa_2\EEE}2|\pp{-}\pp_{\d}|^2+\frac{\SSS\kappa_3\EEE}2|\uu|^2
  \,\d x\d t}
\\[.9em]&\qqqquad\qquad\
\displaystyle{+\int_0^T\!\!\!\int_\varGamma
  \frac{\SSS\varkappa_1\EEE}2|\vv{-}\vv_{\d 1}|^2
  +\frac{\SSS\varkappa_2\EEE}2|\pp{-}\pp_{\d 1}|^2
\,\d S\d t}
\\[.9em]&\qqqquad\qquad\
\displaystyle{+\int_\varOmega\frac{\SSS\lambda_1\EEE}2|\vv(T){-}\vv_{\d T}|^2+\frac{\SSS\lambda_2\EEE}2|\pp(T){-}\pp_{\d T}|^2 \,\d x} &
\\[.6em]
\mbox{subject to}\!\! & \!(\vv,\pp)\text{ satisfying \eq{NS-eqn-evol+}
in the weak sense},
\\[.2em]
& \uu\In 
  {\mathcal U},\ \ \vv\in{\mathcal V},\ \
\pp\in{\mathcal P},
\end{array}
\right.\hspace*{-2em}
\label{NS-eqn-control-semi}\end{align}
\SSS where $\kappa$'s, $\varkappa$'s, and $\lambda$'s are nonnegative constants,
$\kappa_3>0$. \EEE
Here, rather to avoid technicalities, the control $\uu$ is in the bulk, which
makes the problem a bit academical, although some boundary control or control
through initial conditions might be consider\FFF{ed}\EEE, too.

The optimality conditions involves the multipliers
$(\CHI,\LAMBDA)\in{\mathcal V}{\times}{\mathcal P}$ and the Lagrangian
$$
\mathscr{L}(\uu,\vv,\pp,\CHI,\LAMBDA)
:=\large\langle\varPi(\uu,\vv,\pp)\,,(\CHI,\LAMBDA)\large\rangle-
\varPhi(\uu,\vv,\pp)\,.
$$
We denote
$S:\uu\mapsto(\vv,\pp):{\mathcal U}\to{\mathcal V}{\times}{\mathcal P}$
the control-to-state mapping (which was shown in Section~\ref{sec-n=2}
single-valued and continuous) and the composed cost $J(u):=\varPhi(\uu,S(\uu))$.
If $J$ has the G\^ateaux differential $J':{\mathcal U}\to{\mathcal U}^*$,
the first-order optimality conditions for our unconstrained problem
looks simply as $J'(\uu)=0$.

The G\^ateaux differential (if exists) can be calculated by the so-called
adjoint-equation technique. More specifically, 
$J'(\uu)=\mathscr{L}_\uu'(\uu,\vv,\pp,\CHI,\LAMBDA)$ provided
the adjoint velocity $\CHI$ and the adjoint pressure $\LAMBDA$ satisfied
\begin{subequations}\begin{align}\label{OC1}
    \langle\varPi_{(\vv,\pp)}'(\uu,\vv,\pp)\,,(\CHI,\LAMBDA)\rangle
    &=\varPhi_{(\vv,\pp)}'(\uu,\vv,\pp)\,.
    \intertext{Thus, the mentioned first-order
      optimality conditions $J'(\uu)=0$ then result to} 
      \label{OC2}
    \langle\varPi_{\uu}'(\uu,\vv,\pp)\,,(\CHI,\LAMBDA)\rangle
    &=\varPhi_{\uu}'(\uu,\vv,\pp)\,.
\end{align}
\end{subequations}

More specifically, realizing the partial G\^ateaux derivatives
$\varPi_{(\vv,\pp)}'(\uu,\vv,\pp)\in{\rm Lin}({\mathcal V}{\times}{\mathcal P},
      {\mathcal V}^*{\times}{\mathcal P}^*)$
and $\varPhi_{(\vv,\pp)}'(\uu,\vv,\pp)\in
{\rm Lin}({\mathcal V}{\times}{\mathcal P},\R)\cong
      {\mathcal V}^*{\times}{\mathcal P}^*$,
 the adjoint equation \eq{OC1} means
  \begin{align*}
&    \forall (\wt\CHI,\wt\LAMBDA)\!\in\!{\mathcal V}{\times}{\mathcal P}:\ 
\big\langle[\varPi_{(\vv,\pp)}'(\uu,\vv,\pp)](\wt\CHI,\wt\LAMBDA)\,,(\CHI,\LAMBDA)\big\rangle
=\big\langle\varPhi_{(\vv,\pp)}'(\uu,\vv,\pp),(\wt\CHI,\wt\LAMBDA)\big\rangle
  \end{align*}
  which further means
    \begin{subequations}\begin{align}\label{OC1+}
\big[\varPi_{(\vv,\pp)}'(\uu,\vv,\pp)\big]^*(\CHI,\LAMBDA)=\varPhi_{(\vv,\pp)}'(\uu,\vv,\pp)
\intertext{where $[\,\cdot\,]^*$ denotes the adjoint operator, and similarly
\eq{OC2} reads as}
\label{OC2+}
\big[\varPi_{\uu}'(\uu,\vv,\pp)\big]^*(\CHI,\LAMBDA)
    =\varPhi_{\uu}'(\uu,\vv,\pp)\,.
 \end{align}\end{subequations}
    
We note that the control-to-state mapping $S$ is in our two-dimensional case
even continuously differentiable. It is important that the adjoint equation
\eq{OC1+} has a solution, for which it suffices to show that 
$\varPi_{(\vv,\pp)}'(\uu,\cdot,\cdot)$ is surjective
\SSS because then, by open-mapping theorem, there exists a continuous
inverse operator. Then, from the state equation $\varPi(\uu,S(\uu))=0$,
we get $\varPi_{\uu}'(\uu,S(\uu))+\varPi_{(\vv,\pp)}'(\uu,S(\uu))\circ S'(\uu)=0$
so that the G\^ateaux differential of $S$ is given by the explicit formula
$S'(\uu)=[\varPi_{(\vv,\pp)}'(\uu,S(\uu))]^{-1}\varPi_{\uu}'(\uu,S(\uu))$
and it depends continuously on $\uu$. By this surjectivity, also the adjoint
equation \eqref{OC1+} has a solution 
$(\CHI,\LAMBDA)=[\varPi_{(\vv,\pp)}'(\uu,\vv,\pp)]^{-1}(\varPhi_{(\vv,\pp)}'(\uu,\vv,\pp))$.
As it is the only solution of \eqref{OC1+}, the adjoint state is determined
uniquely for a current $\uu$. \EEE

The semi-compressible system contains three bi-linear and one quadratic
terms in \eq{NS-eqn-evol+}, namely
$\varrho(\vv\cdot\!\nablax)\vv$, $\frac\varrho2({\rm div}\,\vv)\,\vv$,
$\frac{\beta}2\nabla(\pp^2)$, and
$\beta \vv\cdot\!\nablax\pp$. These terms give rise \FFF to \EEE seven bilinear terms
in the adjoint system, mixing the state and the adjoint variables.
These seven bilinear terms arise by the Green formula through
the following detailed componentwise calculus
\begin{align}\label{linearized-system}
&\int_\varOmega\bigg(
\CHI{\cdot}\left(\varrho\wt \vv\cdot\!\nablax\vv+\varrho\vv\cdot\!\nablax\wt \vv
+\frac\varrho2({\rm div}\,\vv)\wt \vv+\frac\varrho2({\rm div}\,\wt \vv)\vv
+\beta\nablax(\pp\wt\pp)\right)
\\[-.3em]&\nonumber\hspace*{17em}
+\LAMBDA\big(\wt \vv\cdot\!\nablax\pp+\vv\cdot\!\nablax\wt\pp\big)\bigg)
\,\d x
\linebreak
\\[-.3em]\nonumber&=\sum_{i,j=1}^n\int_\varOmega\!\bigg(\!\CHIj\Big(
\varrho\wt \vv_i\frac{\partial \vv_j\!}{\partial x_i\!}
+\varrho \vv_i\frac{\partial\wt \vv_j\!}{\partial x_i\!}
+\frac\varrho2\frac{\partial \vv_i}{\partial x_i}\wt \vv_j
+\frac\varrho2\frac{\partial\wt \vv_i\!}{\partial x_i\!}\vv_j
+\beta\frac{\partial}{\partial x_j\!}(\pp\wt\pp)
\Big)
\\[-.5em]&\nonumber\hspace*{16.7em}
+\beta\LAMBDA\Big(\wt \vv_i\frac{\partial\pp}{\partial x_i}
+\vv_i\frac{\partial\wt\pp}{\partial x_i}\Big)\bigg)
\,\d x
\\[-.2em]&\nonumber=\sum_{i,j=1}^n\int_\varOmega\!\!\bigg(\!
\Big(\varrho\CHIi\frac{\partial\vv_i}{\partial x_j}
-\varrho\frac{\partial}{\partial x_i\!\!}(\vv_i\CHIj)
+\frac\varrho2\frac{\partial \vv_i}{\partial x_i}\CHIj
-\frac\varrho2\frac{\partial}{\partial x_j\!\!}(\vv_i\CHIi)
+\beta\LAMBDA\frac{\partial\pp}{\partial x_j}
\Big)\wt \vv_j
\\[-.6em]&\nonumber\hspace*{2em}
-\beta\Big(\frac{\partial(\LAMBDA\vv_i)\!\!}{\partial x_i\!\!}
+\pp\frac{\partial\CHIi}{\partial x_i}\Big)\wt\pp\bigg)
\d x\!
+\!\!\sum_{i,j=1}^n\int_\varGamma\varrho \vv_i\CHIj\wt \vv_j\nn_i
+\frac\varrho2\vv_i\CHIi\wt \vv_j\nn_j
+\beta(\LAMBDA\vv_i{+}\pp\CHIi)\wt\pp\nn_i\d S\!\!\!
\\&\nonumber=\!\int_\varOmega\!\Big(\!\big(\varrho(\nablax \vv)^\top\CHI\!
-\varrho{\rm div}(\vv{\otimes}\CHI)
+\frac\varrho2({\rm div}\,\vv)\CHI\!-\frac\varrho2\nablax(\vv{\cdot}\CHI)
+\beta\LAMBDA\nablax\pp\big){\cdot}\wt\vv
\\[-.5em]&\nonumber\hspace*{2em}
-\beta\big({\rm div}(\LAMBDA\vv)+\pp\,{\rm div}\,\CHI\big)\wt\pp\Big)\,\d x
+\int_\varGamma
\Big(\frac{3\varrho}2\CHI{\cdot}\wt \vv+
\beta\wt\pp\LAMBDA\Big)
(\vv{\cdot}\nn)
+\beta\pp\wt\pp(\CHI{\cdot}\nn)
\,\d S\,,
\end{align}
where the boundary integral actually vanishes due to the boundary
conditions $\vv{\cdot}\nn=0$ and $\CHI{\cdot}\nn=0$.
From this, we can read the five bi-linear terms in \eq{QCNS-adjoint-evol1}
when varying the test function $\wt \vv$ and two
bi-linear terms in \eq{QCNS-adjoint-evol2} when varying the test function
$\wt\pp$.
By straightforward modifications of the estimates from Sect.~\ref{sec-n=2}, cf.\
in particular \eq{estimate}--\eq{estimate+}, we can see integrability
of these tri-linear terms in this two-dimensional case.
The resting linear parabolic terms contributing to the linear operator 
$\varPi_{(\vv,\pp)}'(\uu,\cdot,\cdot)$ are standard.
  
More specifically, \eqref{OC1+} leads to the adjoint terminal-boundary-value
for a linear parabolic system
for the multipliers $\CHI$ and $\LAMBDA$ in the classical formulations reads as
\begin{subequations}\label{QCNS-adjoint-evol}\begin{align}\label{QCNS-adjoint-evol1}
    &-\varrho
    \DT\CHI
+\varrho(\nablax \vv)^\top\CHI\!-\varrho
{\rm div}(\vv\otimes\CHI)
+\beta\LAMBDA\nablax\pp-{\rm div}(\nu\EE(\CHI))
&
\\[-.3em]\nonumber
&\hspace*{4em}+\frac\varrho2({\rm div}\,\vv)\CHI-\frac\varrho2\nablax(\vv{\cdot}\CHI)+\nablax\LAMBDA
\,=\,
\SSS\kappa_1(\EEE\vv{-}\vv_{\d})&&\text{in }\ \II{\times}\O,
\nonumber
\\&\label{QCNS-adjoint-evol2}
-\beta\DT\LAMBDA
-\beta{\rm div}(\LAMBDA\vv)
-\GM\Delta\LAMBDA\SSS-\EEE(1{+}\beta\pp){\rm div}\,\CHI
=\,\SSS\kappa_2(\EEE\pp{-}\pp_{\d})
\!\!&&\text{in }\ \II{\times}\O\,,
\\\label{QCNS-adjoint-evol3}
&[\nu\EE(\CHI)\nn]_{\text{\sc t}}^{}+b\vvt
\,=\,\SSS\varkappa_1(\EEE\vv{-}\vv_{\d1})
\ \ \text{ and }\ \ \ 
\nn{\cdot}\CHI=0&&\text{on }\II{\times}\G\,,
\\\label{QCNS-adjoint-evol3+}
&\GM\nablax\LAMBDA{\cdot}\nn=
\SSS\varkappa_2(\EEE\pp{-}\pp_{\d1})&&\text{on }\II{\times}\G\,,
\\\label{QCNS-adjoint-evol4}
&\CHI(T,\cdot)=
\SSS\lambda_1(\EEE\vv(T){-}\vv_{\d T})\,\ \text{ and }\,\ \LAMBDA(T,\cdot)=
\SSS\lambda_2(\EEE\pp(T){-}\pp_{\d T})&&\text{on }\O\,.
\end{align}\end{subequations}
The second condition \eqref{OC2} leads here simply to $\CHI=\SSS\kappa_3\EEE\uu$.

\SSS Since $\kappa_3>0$, the functional $\uu\mapsto\varPhi(\uu,S(\uu))$ is coercive on
${\mathcal U}$. Existence of optimal controls, i.e.\ solutions to
\eq{NS-eqn-control-semi}, can then be shown by the classical direct-method
argument. Here we use the (weak,weak)-continuity of $S$
from Lemma~\ref{lem} and also the (weak,weak)-continuity of the trace operator
$(\vv,\pp)\mapsto(\vv|_{I\times\Gamma}^{},\pp|_{I\times\Gamma}^{}):{\mathcal V}\times{\mathcal P}\to L^2(I{\times}\Gamma;\R^d{\times}\R)$ and also
of the operator $(\vv,\pp)\mapsto(\vv(T),\pp(T)):{\mathcal V}\times{\mathcal P}
\to L^2(\Omega;\R^d{\times}\R)$. The weak lower semicontinuity of
$\uu\mapsto\varPhi(\uu,S(\uu))$ is then implied by convexity of $\varPhi$.

\SSS The existence and uniqueness of the solution to the adjoint problem has been
already discussed on the abstract level, based on the solvability
of the linearized state equation for all right-hand sides, i.e.\ surjectivity
of the G\^ateaux derivative of the state equation with respect to $(\vv,\pp)$.
At a given $(\vv,\pp)$, this linearization in the direction $(\wt\vv,\wt\pp)$
gives seven terms originated from the bilinear/quadratic nonlinearities in
(\ref{NS-eqn-evol+}a,b). The a-priori estimates for this linearized system
can be performed by a test by $(\wt\vv,\wt\pp)$, which gives the trilinear
terms of the type ${\rm div}\vv|\wt\vv|^2$, $(\vv{\cdot}\nablax)\wt\vv{\cdot}\wt\vv$,
$(\vv{\cdot}\nablax\pp)\pp$, and $({\rm div}\vv)\wt\pp^2$,
cf.\ \eqref{linearized-system} for $\CHI=\wt\vv$ and $\LAMBDA=\wt p$.
These terms can be estimated by Gagliardo-Nirenberg, Korn, H\"older, Young,
and Gronwall inequalities quite similarly as we did in \eq{unique-est}. Again, it
works only for the two-dimensional situation.
\EEE

Let us briefly summarize the above arguments and calculations:

\newpage

\begin{proposition}[Existence and optimality conditions]\label{prop1}
Let $\varrho,\nu,\beta,\gamma>0$, \eq{IC} hold, and $n=2$. Then:\\
    \ITEM{(i)}{The optimal-control problem \eq{NS-eqn-control-semi} possesses
     at least one solution $(\uu,\vv,\pp)$.}
    \ITEM{(ii)}{For any such solution $(\uu,\vv,\pp)$, there exists the \SSS
      unique \EEE adjoint
      state $(\CHI,\LAMBDA)\in{\mathcal V}{\times}{\mathcal P}$ satisfying
      the terminal-boundary-value problem \eq{QCNS-adjoint-evol} in the
      weak sense and $\uu=\CHI\SSS/\kappa_3\EEE$.
In particular, any optimal control $\uu\in{\mathcal U}$ must be
more regular, belonging also to ${\mathcal V}$ from \eq{spaces-V}.}
\end{proposition}

\section{Concluding remarks}\label{sec-remarks}

We end this article by several remarks suggesting modifications
or enhancements of the controlled semilinear parabolic system
still keeping the bilinear character of all involved nonlinearities
and widening applicability towards three-dimensional situations
or coupling with other phenomena.

\begin{remark}[{\sl Multipolar fluids for three-dimensional problems}]\label{rem-nonsimple}
\upshape
The semi-compressible model (similarly as fully compressible for $\beta=0$ and
$\gamma=0$) works only in two-dimensional situations because the needed uniqueness of
the response is not granted in higher dimensions. To extend the above results 
towards 3-dimensional problems, a concept of so-called 
multipolar fluids can be used. This introduces a higher-order
viscosity. Under the name 2nd-grade nonsimple fluids,
\eq{NS-eqn-evol-mult} was devised by E.\,Fried and M.\,Gurtin \cite{FriGur06TBBC}
and earlier, even more generally and nonlinearly as multipolar fluids, by
J.\,Ne\v cas at al.\ \cite{BeBlNe92PBMV,NeNoSi89GSIC,NecRuz92GSIV}.
Here it means that (\ref{NS-eqn-evol+}a,b) modifies as
\begin{subequations}\label{NS-eqn-evol-mult}\begin{align}\label{NS-eqn-evol1-mult}
&\varrho\DT\vv+\varrho(\vv\cdot\!\nablax)\vv
    -{\rm div}\big(\nu\EE(\vv)-{\rm div}(\nu_1\nablax\EE(\vv))\big)
+\nablax\Big(\pp+\frac\beta2\pp^2\Big)=\uu
\hspace*{-12em}&&
\\[-.4em]\nonumber&\hspace*{7em}\text{ and }\ \ \ \ \ 
\beta(\DT\pp+\vv\cdot\!\nablax\pp)
+{\rm div}\,\vv=\GM\Delta\pp&&\text{on }\II{\times}\O\,,
\\[-.2em]\label{NS-eqn-evol2-mult}
&\big[\nu\EE(\vv)\nn
  -\divS(\nu_1\nablax\EE(\vv))\nn\big]_{\text{\sc t}}^{}\!
+b\vvt
\,=\,0\,,\ \ \ \ \ 
\nn{\cdot}\vv=\,0\hspace*{-2em}
\\[-.6em]&\hspace*{5em}\text{ and }\ \ \ 
\nablaxx \vv\,\vdots\,(\nn{\otimes}\nn)=0,\ \ \ 
 \nn{\cdot}\nablax\pp=0
&&\text{on }\II{\times}\G,
\nonumber\end{align}\end{subequations}
\TTT with a (presumably small) ``hyper''-viscosity coefficient $\nu_1>0$ and \EEE
with the initial condition\TTT{s} \EEE \eq{NS4+}. In \eq{NS-eqn-evol2-mult},
``$\divS$'' denotes the surface divergence defined as
$\divS(\cdot)=\mathrm{tr}(\nablaS(\cdot))$ with
$\mathrm{tr}(\cdot)$ denoting the trace and $\nablaS$ denoting the
surface gradient given by $\nablaS v=(\mathbb I- \nn{\otimes}\nn)\nabla v
=\nabla v-\frac{\partial v}{\partial\nn}\nn$.
The energetics \eq{energy-NS-semi} expands by the dissipation rate
$\nu_1|\nablaxx \vv|^2$ so that $\mathcal V$ from \eq{spaces-V} takes
$H^2(\O;\Rn)$ instead of $H^1(\O;\Rn)$.
The above uniqueness arguments can now be modified for $n=3$
too, cf.\ \cite[Prop.\,4]{Roub20QISI}, and again $\varPi(\uu,\cdot,\cdot):
{\mathcal V}{\times}{\mathcal P}\to{\mathcal V}^*{\times}{\mathcal P}^*$.
Here one uses that, by interpolation,
$L^\infty(I;L^2(\Omega))\cap L^2(I;H^2(\Omega))\subset L^4(I;H^1(\Omega))
\subset L^4(I;L^6(\Omega))$ so
that \eq{estimate} holds in the modification
\begin{align}\nonumber
\|(\vv{\cdot}\!\nablax)\vv{\cdot}\wt\vv\|_{L^1(I{\times}\O)}^{}\!\le
C\|\vv\|_{L^4(I;L^6(\O;\R^3))}^{}\|\nablax\vv\|_{L^4(I;L^2(\O;\R^{3\times3}))}^{}
\|\wt\vv\|_{L^4(I;L^6(\O;\R^3))}^{}
\end{align}
while the nonlinear term $\pp^2{\rm div}\,\wt\vv$ bears the estimation
\begin{align}\nonumber
  &\|\pp^2{\rm div}\,\wt\vv\|_{L^1(I{\times}\O)}^{}\le
  \|\pp\|_{L^{8/3}(I;L^4(\O))}^2\|{\rm div}\,\wt\vv\|_{L^4(I;L^2(\O))}^{}
  \\&\nonumber\qquad\le C
  \|\pp\|_{L^\infty(I;L^2(\O))}^{1/2}\|\pp\|_{L^2(I;H^1(\O))}^{3/2}
  \|\wt\vv\|_{L^\infty(I;L^2(\O;\R^3))}^{1/2}\|\wt\vv\|_{L^2(I;H^2(\O;\R^3))}^{1/2}
\end{align}
when using twice the Gagliardo-Nirenberg interpolation;
note that here the estimation is exact without any ``reserve''. Similar
estimate hold for the nonlinear term $(\vv\cdot\!\nabla\pp)\wt\pp$,
namely
\begin{align}\nonumber
  \|(\vv\cdot\!\nabla\pp)\wt\pp\|_{L^1(I{\times}\O)}^{}&\le
  \|\vv\|_{L^4(I;L^6(\O;\R^3))}^{}\|\nabla\pp\|_{L^2(I{\times}\O)}^{}\|\wt\pp\|_{L^4(I;L^3(\O))}^{}
  \\\nonumber& \le C\|\vv\|_{L^\infty(I;L^2(\O;\R^3))}^{1/2}\|\vv\|_{L^2(I;H^2(\O;\R^3))}^{1/2}\times
  \\&\nonumber\hspace*{1em}\times\|\nabla\pp\|_{L^2(I{\times}\O)}^{}
  \|\wt\pp\|_{L^\infty(I;L^2(\O))}^{1/2}\|\wt\pp\|_{L^2(I;L^6(\O))}^{1/2}
\end{align}
due to the Gagliardo-Nirenberg interpolation; note that again the estimation is exact
without any ``reserve''. Thus Proposition~\ref{prop1} holds for $n=3$
for this 2nd-grade nonsimple semi-compressible system \eqref{NS-eqn-evol-mult}.
\end{remark}

\begin{remark}[{\sl Multipolar fluids for pressure constraints}]
  \upshape\label{NS-mult3}
  One can think about 3rd-grade nonsimple semi-compressible fluids governed
  by the system:
  \begin{subequations}\begin{align}
&\varrho\DT\vv+\varrho(\vv\cdot\!\nablax)\vv
-{\rm div}\big(\nu\EE(\vv){+}{\rm div}^2(\nu_1\nabla^2\EE(\vv))\big)
+\nablax\Big(\pp+\frac\beta2\pp^2\Big)=\,\uu\,,
\\[-.2em]&\label{NS-mult3-p}\beta(\DT\pp
+\vv\cdot\!\nablax\pp)+{\rm div}\,\vv=0 
\end{align}\end{subequations}
on $I{\times\O}$ with the corresponding three boundary conditions for $\vv$ on
$I{\times}\G$; these conditions
are rather technical and we refer to \cite{Roub20CHEC}. On the other hand,
a reasonable analysis can be performed without pressure diffusion $\GM\Delta\pp$
in \eqref{NS-mult3-p} and thus no boundary condition is prescribed for pressure.
Now $\vv\In C_{\rm w}(I;H^1(\varOmega;\R^n))\:\cap\:L^2(I;H^3(\varOmega;\R^n))$, so
that even for $n=3$, we now have
$\nabla\vv\in L^2(I;L^\infty(\varOmega;\R^{3\times3}))$.
Then, if the initial pressure is enough regular, this regularity now will
be transported along the evolution. Namely, if $\pp_0\in W^{1,q}(\varOmega)$,
we obtain $\pp\in L^\infty(I;W^{1,q}(\varOmega))$ by applying $\nabla$ to 
\eqref{NS-mult3-p} and then testing it by $|\nabla\pp|^{q-2}\nabla\pp$. Here
we use the calculus
\begin{align}\nonumber
&\frac1q\frac{\d}{\d t}\int_\varOmega|\nabla\pp|^q\,\d x
=\int_\varOmega|\nabla\pp|^{q-2}\nabla\pp\cdot\nabla\DT\pp\,\d x
\\&\nonumber
\qquad
=-\int_\varOmega|\nabla\pp|^{q-2}\nabla\pp{\cdot}\big(\nabla(\vv\cdot\!\nablax\pp)
+\nabla{\rm div}\,\vv\big)\,\d x
\\&\nonumber
\qquad=\int_\varOmega\EE(\vv){:}\Big(|\nabla\pp|^{q-2}\nabla\pp{\otimes}\nabla\pp
-\frac1q|\nabla\pp|^q\mathbb I\Big)-|\nabla\pp|^{q-2}\nabla\pp\cdot\!\nabla{\rm div}\,\vv\,\d x
\\&\nonumber
\qquad\le2\|\EE(\vv)\|_{L^\infty(\varOmega;\R^{n\times n})}^{}\|\nablax\pp\|_{L^q(\varOmega;\R^n)}^q\!
\\&\hspace{13em}
+\|\nabla\EE(\vv)\|_{L^{q}(\varOmega;\R^{n\times n\times n})}^{q}\!
+\|\nabla\pp\|_{L^{q}(\varOmega;\R^n)}^{q}
\nonumber\end{align}
and then Gronwall's inequality. The penultimate term needs $q\le 2n/(n{-}2)$
or $q<\infty$ for $n=2$. Also, one can see the estimate
$\DT\pp\in L^2(I{\times}\O)$.
Notably, for $q>n$, we have pressure $\pp$ in $C(I{\times}\bar\varOmega)$ and,
in particular, also traces on $\varGamma$ are in $C(I{\times}\varGamma)$.
Therefore, one can consider also the pointwise state constraints on pressure
even on the boundary, e.g.\ a technologically relevant condition on the
local pressure on the wall $\varGamma$ of the container $\varOmega$ of the type
$|\pp|_\varGamma|\le\pp_{\rm max}$. This set of $\pp$'s has a nonempty interior in
$C(I{\times}\varGamma)$ and thus the corresponding multiplier in the 1st-order
condition is well determined as a measure on $I{\times}\varGamma$.
\end{remark}

\begin{remark}[{\sl Enhancement: Cahn-Hilliard diffusion}]
  \upshape
  Liquids can contain some other constituent (e.g.\ salt in
  sea water or Nickel in molten Iron outer Earth core).   
  This consistent with a concentration $\chi$ can diffuse according 
  the gradient of a chemical potential $\mu$. Then the
  original semi-compressible initial-boundary-value system
  \eqref{NS-eqn-evol+} expands as:
 \begin{subequations}\label{CH-eqn-evol+}\begin{align}\label{CH1+}
&\varrho\DT\vv
     +\varrho(\vv\cdot\!\nablax)\vv-{\rm div}\big(\nu\EE(\vv)
     +\alpha\nablax\chi{\otimes}\nablax\chi\big)
\\[-.3em]&\nonumber\hspace*{0em}
+\frac\varrho2({\rm div}\,\vv)\,\vv+\nablax\Big(\pp+\frac\beta2\pp^2\!+
(\chi{-}\chi_{\rm eq})^2\!+\frac\alpha2|\nablax\chi|^2\Big)
\,=\,\uu\!\!\!\!
&&\text{on }\II{\times}\O\,,
\\[-.2em]\label{CH2+}
&\beta(\DT\pp+\vv\cdot\!\nablax\pp)
+{\rm div}\,\vv=\GM\Delta\pp&&\text{on }\II{\times}\O\,,
\\[-.2em]\label{CH3+}
&\DT\chi+\vv\cdot\!\nablax\chi={\rm div}(m\nablax\mu)\ \ \text{ with }\ \ 
\mu=(\chi{-}\chi_{\rm eq})^2-\alpha\Delta\chi
&&\text{on }\II{\times}\O\,,
\\[-.2em]&\big[\nu\EE(\vv)\nn\big]_{\text{\sc t}}^{}\!+b\vvt=0,\ \ \  
\nn{\cdot}\vv=0,\ \ \nn{\cdot}\nablax\pp=0,\ \ \nn{\cdot}\nablax\chi=0\!
&&\text{on }\II{\times}\G\,,
\label{CH4+}
\\[-.2em]
&\vv(0,\cdot)=\vv_0\,,\ \ \ \pp(0,\cdot)=\pp_0\,,\ \ \text{ and }\ \ \ \chi(0,\cdot)=\chi_0 &&\text{on }\O\,,
 \end{align}\end{subequations}
where $\chi_{\rm eq}$ is an equilibrium concentration, 
 $m>0$ a mobility and $\alpha>0$ a capillarity coefficient.
 The equation \eqref{CH3+} is called the Cahn-Hilliard one, and 
 the extra symmetric stress $\alpha\nablax\chi{\otimes}\nablax\chi-
 \frac\alpha2|\nablax\chi|^2\mathbb I$ occurring in \eqref{CH1+}
 is called Korteweg's stress.
 The energy balance \eq{energy-NS-semi} now augments as 
 \begin{align}\label{energy-NS-semi}
   &\int_\varOmega\!\!\!\!\linesunder{\frac\varrho2|\vv(t)|^2_{_{_{_{_{_{_{}}}}}}}}{kinetic}{energy}\!\!\!\!\!\!+\!\!\!\!\!\!\linesunder{\frac{\beta}2\pp(t)^2+
     (\chi{-}\chi_{\rm eq})^2
+\frac\alpha2|\nablax\chi|^2
     _{_{_{_{_{_{_{}}}}}}}}{stored}{energy}\!\!\!\!\!\d x
\\[-.1em]&\nonumber\ \ \ \ 
+\int_0^t\!\int_\varOmega\!\!\!\linesunder{\nu|\EE(\vv)|^2+\GM|\nablax\pp|^2+m|\nabla\mu|^2_{_{_{_{_{_{_{}}}}}}}}{dissipation rate}{in the bulk}\!\!\!\d x\d t
  +\int_0^t\!\int_\varGamma\!\!\!\!\!\!\!\!\!\!\!\linesunder{
 b|\vvt|^2{}_{_{_{_{_{_{_{}}}}}}}}{dissipation rate}{on the boundary}\!\!\!\!\!\!\!\!\!\!\!\!\!
\d S\d t 
\\[-.1em]&\nonumber\hspace*{0em}=\int_0^t\!\int_\varOmega\!\!\!\!\!\!\!\linesunder{
  \uu{\cdot}\vv_{_{_{_{_{_{_{}}}}}}}}{power
  of}{the control}\!\!\!\!\!\!\!\!\!\!\d x\d t
  +\int_\varOmega\frac\varrho2|\vv_0|^2+\frac{\beta}2\pp_0^2
   +(\chi_0{-}\chi_{\rm eq})^2
+\frac\alpha2|\nablax\chi_0|^2\,\d x.
\nonumber\end{align}
From this, a similar analysis of the controlled system and an optimal
control problem can be casted similarly as in Sections~\ref{sec-system}
and \ref{sec-control}. Optimal control for the incompressible variant in two
dimensions has been treated in \cite{Medj15OCCH} and in a certain 
nonlocal variant \cite{BiDhMo20MPSO,FrRoSp16ODCN}.
\end{remark}

\begin{remark}[{\sl Enhancement: magneto-hydrodynamics}]
  \upshape
  Some fluids are electrically conductive. Typically it concerns
  molten metals, like hot Iron with Nickel in the outer core of Earth
  or metallic hydrogen in Jupiter and Saturn. It maybe also plasma
  especially in stellar applications, which
  was the original motivation for this model. This brings an interesting
  coupling of semi-compressible fluids
  with the magnetic induction $\bb$. The other parameter is 
  electric conductivity $\sigma$. The so-called induction equation
  merges Faraday's law and Ohm's law: 
\begin{align*}
  \DT\bb={\rm rot}(\vv{\times}\bb)+
{\rm rot}\,\Big(\frac1{\mu_0\sigma}{\rm rot}\,\bb\Big)
  \ \ \text{ and }\ \ {\rm div}\,\bb=0\,,
\end{align*}
where $\mu_0$ is the vacuum permeability or, using the calculus
${\rm rot}(\vv{\times}\bb)=(\bb{\cdot}\nabla)\vv-(\vv\cdot\!\nabla)\bb$, also
$$
\DT\bb+(\vv\cdot\!\nablax)\bb-(\bb\,{\cdot}\nabla)\vv=
{\rm rot}\,\Big(\frac1{\mu_0\sigma}{\rm rot}\,\bb\Big)
\ \ \ \text{ and }\ \ {\rm div}\,\bb=0\,.
$$
When $\sigma$ is constant,
then we can further use the calculus ${\rm rot}\,{\rm rot}\,\bb
=\nabla({\rm div}\,\bb)-\Delta\bb=-\Delta\bb$, so that
\begin{align*}
\DT\bb+(\vv\cdot\!\nablax)\bb
-(\bb\,{\cdot}\nabla)\vv=
\frac1{\mu_0\sigma}\Delta\bb\ \ \text{ and }\ \ {\rm div}\,\bb=0\,.
\end{align*}
The magnetic field influences the mechanical part through
Lorenz' force which, under electroneutrality, is ${\bm f}={\bm j}{\times}\bb$
with low-frequency Amp\'ere's law neglecting the displacement current so that
$\mu_0{\bm j}={\rm rot}\,\bb$.
Using the calculus $\frac12\nabla(\bb\cdot\bb)=(\bb\cdot\nabla)\bb+\bb{\times}({\rm rot}\,\bb)$, we eventually have
$$
{\bm f}={\bm j}{\times}\bb=
\frac{(\bb\,{\cdot}\nabla)\bb}{\mu_0}-\nabla\frac{|\bb|^2}{2\mu_0}\,.
$$
The semi-compressible system \eq{NS-eqn-evol+} now expands as:
\begin{subequations}\label{MHD-eqn-evol+}
\begin{align}\label{MHD1+}
&\varrho\DT\vv
     +\varrho(\vv\cdot\!\nablax)\vv-{\rm div}\big(\nu\EE(\vv)\big)
\\[-.3em]&\nonumber\hspace*{1em}
+\frac\varrho2({\rm div}\,\vv)\,\vv+
\nablax\Big(\pp+\frac\beta2\pp^2\!+\frac{|\bb|^2}{2\mu_0}
\Big)
\,=\,\frac{(\bb\,{\cdot}\nablax)\bb}{\mu_0}+\uu\!\!\!
&&\text{on }\II{\times}\O\,,
\\[-.2em]\label{MHD2+}
&\beta(\DT\pp+\vv\cdot\!\nablax\pp)
+{\rm div}\,\vv=\GM\Delta\pp&&\text{on }\II{\times}\O\,,
\\[-.2em]\label{MHD3+}
&\DT\bb+(\vv\cdot\!\nablax)\bb-(\bb\,{\cdot}\nabla)\vv=
  \frac1{\mu_0\sigma}\Delta\bb\ \ \ \text{ and }\ \ {\rm div}\,\bb=0
&&\text{on }\II{\times}\O\,,
\\[-.2em]
&\big[\nu\EE(\vv)\nn\big]_{\text{\sc t}}^{}\!+b\vvt=0,\ \ \  
\nn{\cdot}\vv=0,\ \ \nn{\cdot}\nablax\pp=0,\ \ \nn{\cdot}\nablax\bb=0
&&\text{on }\II{\times}\G\,,
\label{MHD4+}
\\[-.2em]
&\vv(0,\cdot)=\vv_0\,,\ \ \ \pp(0,\cdot)=\pp_0\,,\ \ \text{ and }\ \ \
\bb(0,\cdot)=\bb_0 &&\text{on }\O\,.
\end{align}\end{subequations}
The magneto-hydrodynamic system is the basic model for a magnetic dynamo effect,
and its usage in planetary physics explains the magnetic
field generation in particular in our planet Earth. Mostly this system is
considered incompressible but sometimes the semi-compressible variant
(but without the pressure diffusion term $\gamma\Delta\pp$) can be found in
literature, too, viz \cite{Bisk93NMHD,Schn09LMHD}.
Again, the departure point for analysis is the energy balance like 
\eq{energy-NS-semi} which now involves also the induction
equation tested by the intensity of magnetic field $\hh=\bb/\mu_0$,
so it augments as
 \begin{align}\nonumber
   &\int_\varOmega\!\!\!\!\linesunder{\frac\varrho2|\vv(t)|^2_{_{_{_{_{_{_{}}}}}}}}{kinetic}{energy}\!\!\!+\!\!\!\!\!\!\!\!\!\linesunder{\frac{\beta}2\pp(t)^2
     \ +\ \frac{|\bb|^2}{2\mu_0}
   }{stored energy by pressure}{and magnetic field}\!\!\!\!\!\d x
\\[-.1em]&\nonumber\ \ \ \ 
+\int_0^t\!\int_\varOmega\!\!\!\linesunder{\nu|\EE(\vv)|^2+\GM|\nablax\pp|^2
  +\frac{|\nabla\bb|^2}{\mu_0\sigma}
}{dissipation rate}{in the bulk}\!\!\!\d x\d t
  +\int_0^t\!\int_\varGamma
  \!\!\!\!\!\!\!\!\!\!\!\linesunder{
    b|\vvt|^2{}_{_{_{_{_{_{_{}}}}}}}}{dissipation rate}{on the boundary}\!\!\!\!\!\!\!\!\!\!\!\!\!
\d S\d t 
\\[-.1em]&\nonumber\hspace*{7em}=\int_0^t\!\int_\varOmega\!\!\!\!\!\!\!\linesunder{
  \uu{\cdot}\vv_{_{_{_{_{_{_{}}}}}}}}{power
  of}{the control}\!\!\!\!\!\!\!\!\!\!\d x\d t
  +\int_\varOmega\frac\varrho2|\vv_0|^2+\frac{\beta}2\pp_0^2
  +\frac{|\bb_0|^2}{2\mu_0}\,\d x.
\nonumber\end{align}
From this, we now also read the apriori bound for
$\bb\in L^\infty(I;L^2(\varOmega;\R^n))\,\cap\,L^2(I;H^1(\varOmega;\R^n))$.
The convergence of an (unspecified) approximate solutions in bi-linear terms
$(\vv\cdot\!\nabla)\,\vv$, $({\rm div}\,\vv)\vv$,
$(\vv\cdot\!\nabla)\,\pi$, $\nabla\pi^2$,  $\nabla|\bb|^2$,
$(\bb\,{\cdot}\nabla)\bb$, $(\vv\cdot\!\nabla)\bb$, and
$(\vv{\,\cdot}\nabla)\bb$ is then easy by
the Aubin-Lions compact-embedding theorem.
In the nonsimple variant as in Remarks~\ref{rem-nonsimple} and \ref{NS-mult3},
even we have a rigorous energy conservation because $\DT\bb\in L^2(I;H^1(\varOmega;\R^n)^*)+L^1(I;L^2(\O;\Rn))$ is in duality with
$\bb\in L^2(I;H^1(\varOmega;\R^n))\cap L^\infty(I;L^2(\O;\Rn))$.
The application of the optimization is in optimal control of plasma in
tokamaks or in identification of existing flows in the planetary or stellar
astrophysics.
\end{remark}

\bigskip\bigskip

\baselineskip=10pt
{\small
  \noindent{\it Acknowledgments.}
 The author acknowledges a partial support of the
CSF (Czech Science Foundation) project 19-04956S,
the M\v SMT \v CR (Ministery of Education of the Czech Rep.) project
CZ.02.1.01/0.0/0.0/15$\underline{\ }$003/0000493,
and the institutional support RVO: 61388998 (\v CR).
}

\bigskip\bigskip
\bigskip
\end{marginpar}
\end{document}